\documentclass[a4paper,pdftex,reqno,10pt,english]{amsart}
\usepackage{babel}
\usepackage{epsfig,graphicx,times}
\usepackage{amssymb,amsmath}

\def\Ps{\mathcal{P}}

\def\diam{\mbox{diam}}
\def\car{\mbox{char}}
\def\pfo{\downarrow\!\!}
\newcommand{\Tsm}{\hspace*{0.6cm}}
\newcommand{\di}{m}

\def\@begintheorem#1#2{\list{}{\thm@body}%
  \item[]{\bf #1~#2.}\quad\it\ignorespaces}
\def\@opargbegintheorem#1#2#3{\list{}{\thm@body}%
  \item[]{\bf #1~#2~\ifrembrks #3\global\rembrksfalse\else (#3)\fi.}%
  \quad\it\ignorespaces}
\def\@endtheorem{\endlist}
\newtheorem{theorem}{Theorem}
\newtheorem{lemma}{Lemma}
\newtheorem{definition}{Definition}

\newtheorem{algo}{Algorithm}
\newtheorem{conjecture}{Conjecture}

\begin{document}

\title[On the minimum diameter of plane integral point sets]{On the minimum diameter of plane integral point sets}

\author{Sascha Kurz}
\address{Sascha Kurz\\Department of Mathematics, Physic and Informatics\\University of Bayreuth\\Germany}
\email{sascha.kurz@uni-bayreuth.de}

\author{Alfred Wassermann}
\address{Alfred Wassermann\\Department of Mathematics, Physic and Informatics\\University of Bayreuth\\Germany}
\email{alfred.wassermann@uni-bayreuth.de}

\begin{abstract}
    Since ancient times mathematicians consider geometrical objects with integral side lengths. We consider
    plane integral point sets $\Ps$, which are sets of $n$ points in the plane with pairwise 
    integral distances where not all the points are collinear.

    The largest occurring distance is called its diameter. Naturally the question about the minimum possible 
    diameter $d(2,n)$ of a plane integral point set consisting of $n$ points arises. We give some new exact values 
    and describe state-of-the-art algorithms to obtain them. It turns out that plane integral point sets with
    minimum diameter consist very likely of subsets with many collinear points. For this special kind of point sets
    we prove a lower bound for $d(2,n)$ achieving the known upper bound $n^{c_2\log\log n}$ up to a constant in 
    the exponent.

    A famous question of Erd\H{o}s asks for plane integral point sets with no $3$ points on a line and no 
    $4$ points on a circle. Here, we talk of point sets in general position and denote the corresponding 
    minimum diameter by $\dot{d}(2,n)$. Recently $\dot{d}(2,7)=22\,270$ could be determined via an exhaustive search. 
\end{abstract}

\keywords{integral distances, diameter, exhaustive search, orderly generation}
\subjclass[2000]{52C10; 11D99, 53C65}

\maketitle

\section{Introduction}
  In radio astronomy systems of antennas are used. To avoid frequency losses the distance between each pair of 
  antennas have to be an integer multiple of the used wave length. So there is some interest in the construction 
  and properties of $m$-dimensional integral point sets $\Ps$, i.e. sets of $n$ points in the Euclidean space 
  $\mathbb{E}^m$ with pairwise integral distances where not all the points are contained in a hyperplane of 
  $\mathbb{E}^m$. For other applications we refer to \cite{integral_distances_in_point_sets}.

  In this article we focus on the planar case $m=2$ and refer to \cite{phd_kurz,characteristic,paper_laue} for $m\ge 3$. 
  At most $n-1$ points are allowed to be collinear. A point set is said to be in semi-general position if no three points 
  are collinear. If additionally no four points are located on a circle we talk of general position.
  To describe integral point sets $\Ps$ we denote the largest occurring distance 
  as its diameter $\diam(\Ps)$. From the combinatorial point of view there is a natural interest in the determination 
  of the minimum possible diameter $d(2,n)$ of plane integral point sets. For plane integral point sets in semi-general 
  position and general position we denote the minimum diameter by $\overline{d}(2,n)$ and $\dot{d}(2,n)$, respectively.

  Although the study of integral point sets has a long history, see \cite{integral_distances_in_point_sets} for
  an overview, not much has been known about the exact values of $d(2,n)$, $\overline{d}(2,n)$, and $\dot{d}(2,n)$,
  previously. The first exact values are \cite{integral_distances_in_point_sets,hab_kemnitz}:

  $$\left(d(2,n)\right)_{n=3,\dots ,9}=1,4,7,8,17,21,29,$$ 
  $$\left(\overline{d}(2,n)\right)_{n=3,\dots ,9}=1,4,8,8,33,56,56,$$ 
  and
  $$\left(\dot{d}(2,n)\right)_{n=3,\dots ,6}=1,8,73,174.$$ 

  Apart from $d(2,n)\le\overline{d}(2,n)\le\dot{d}(2,n)$ the best known bounds are given in
  \cite{note_on_integral_distances} and \cite{minimum_diameter}, respectively,
  $$c_1n\le d(2,n)\le\overline{d}(2,n)\le n^{c_2\log\log n}\,.$$

  For pictures of the corresponding point sets we refer to \cite{integral_distances_in_point_sets,kreisel,phd_kurz}.
  It is worth noting that the points of all known integral points sets in semi-general position with minimum diameter
  lie on a common circle each. We remark that there are constructions known which lead to a dense subset of points of 
  the unit circle with pairwise rational distances, see i.e. \cite{integral_distances_in_point_sets}.
  A famous Erd\H{o}s problem asks for seven points in the plane no three on a line, no four on a circle
  with pairwise integral distances or more generally for the existence and value of $\dot{d}(2,n)$ (still open).
%  For $n=3$ the answer is the equilateral triangle with side lengths $1$ and for $n=4,5,6$ the corresponding point 
%  sets are depicted in Figure \ref{fig_n4oc_4_5_6}. No  example for seven points is known so far. Via exhaustive
%  enumeration Kemnitz \cite{hab_kemnitz} has obtained the bound $\dot{d}(2,7)>320$. On the contrary 
  There are some constructions for infinite series of plane integral point sets in general position with $6$ points 
  and coprime side lengths \cite{hab_kemnitz}. Recently such an example consisting of seven points was found
  \cite{kreisel}.

  This is the story of plane integral point sets so far. The paper is arranged as follows:

  In Section \ref{sec_exhaustive_search} we give algorithms for the exhaustive generation of plane integral 
  point sets. A variant of orderly generation which combines two integral point sets to obtain a new one instead of
  extending one point set is described in Subsection \ref{sub_sec_orderly}. Clique search is utilized in Subsection
  \ref{sub_sec_clique_search}. In Section
  \ref{section_characteristic} we analyze properties of the characteristic $\car(\Ps)$, this is the squarefree
  part of $(a+b+c)(a+b-c)(a-b+c)(-a+b+c)$ of a non-degenerated subtriangle with side lengths $a$, $b$, and $c$ of a
  plane integral point set $\Ps$. Subsections \ref{subsection_general_position}, \ref{subsection_semi_general_position},
  and \ref{sub_sec_arbitrary} are devoted to computational results using the algorithms described Section
  \ref{sec_exhaustive_search}:

  \begin{theorem}
    \label{thm_n3ol}
    For $n\le 36$ the plane integral point sets in semi-general position with $n$ points and minimum diameter 
    consist of points on a circle with radius $r=\frac{z}{k}\sqrt{k}$ where $k\in\{3,7,15\}$ and $z$ is an 
    integer with many prime factors.
  \end{theorem}

  \begin{conjecture}
    \label{conjecture_on_circle}
    The points of a plane integral point set in semi-general position are situated on a circle.
  \end{conjecture}

  \begin{theorem}
    \label{thm_general_position}
    $$\dot{d}(2,7)=22\,270.$$
  \end{theorem}

  This improves the previous known bound $\dot{d}(2,7)>320$ \cite{hab_kemnitz} and answers Erd\H{o}s' question 
  positively. The corresponding integral point set was first announced in \cite{kreisel}.

%  In Subsection \ref{sub_sec_clique_search} we reformulate the construction of plane integral point sets with minimum 
%  diameter as a clique search problem. The computational results are given in Section \ref{sec_arbitrary}, where 
%  we determine the exact values of $d(2,n)$ up to $n=122$.

  \begin{theorem}
    \label{thm_on_line}
    For $9\le n\le 122$ the plane integral point sets with $n$ points and minimum diameter $d(2,n)$ consist of a subset 
    of $n-1$ collinear points and one point apart from this line.
  \end{theorem}

  \begin{conjecture}
    \label{conjecture_on_line}
    For $n\ge 9$ a plane integral point set with minimum diameter contains a subset of $n-1$ collinear points
  \end{conjecture}

  Theorem \ref{thm_on_line} motivates us to investigate plane integral point sets with many collinear points in Section
  \ref{sec_decomposition_lemma}. We give a link between these special plane integral point sets and factorizations of 
  integers. This link enables us to generate them very efficiently and to give the following lower bound.

  \begin{theorem}
    \label{thm_lower_bound}
    For $\delta>0$, $\varepsilon>0$, and a plane integral point set $\mathcal{P}$ consisting of $n$ points with at
    least $n^\delta$ collinear points there exists a $n_0(\varepsilon)$ such that for all $n\ge n_0(\varepsilon)$ we have
    $$\mbox{diam}(\mathcal{P})\ge n^{\frac{\delta}{4\log 2(1+\varepsilon)}  \log\log n}.$$
  \end{theorem}

We end with some remarks on lower bounds for $d(2,n)$ in Section \ref{sec_some_remarks} and give a conclusion.

%  The results lead to a strategy to prove a lower bound $d(2,n)>n^{c\log\log n}$ in the arbitrary case described 
%  in Section \ref{sec_proof_strategy}.
%
% ---------------------------------------------------------------------------------------------------------------------  
\section{Exhaustive search}
\label{sec_exhaustive_search}

\noindent
To determine some further exact values of $d(2,n)$, $\overline{d}(2,n)$, and $\dot{d}(2,n)$ we have applied an exhaustive search. In this Section we describe the used algorithms.

\subsection{Orderly generation by combination}
\label{sub_sec_orderly}

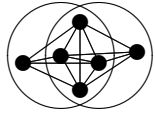
\begin{figure}[ht]
    \begin{center}
      \setlength{\unitlength}{0.5cm}
      \begin{picture}(5,2)
        \put(2,.5){\circle{4}}
        \put(3,.5){\circle{4}}
        \put(2,.5){\circle*{0.4}}
        \put(2.5,-0.4){\circle*{0.4}}
        \put(2.5,1.4){\circle*{0.4}}
        \put(3,0.3){\circle*{0.4}}
        \put(1,0.3){\circle*{0.4}}
        \put(4,0.6){\circle*{0.4}}
        \put(0,-1.5){\qbezier(2,2)(2.25,1.55)(2.5,1.1)}
        \put(0,-1.5){\qbezier(2,2)(2.25,2.45)(2.5,2.9)}
        \put(0,-1.5){\qbezier(2,2)(2.5,1.9)(3,1.8)}
        \put(0,-1.5){\qbezier(2,2)(1.5,1.9)(1,1.8)}
        \put(0,-1.5){\qbezier(2.5,1.1)(2.5,2.0)(2.5,2.9)}
        \put(0,-1.5){\qbezier(2.5,1.1)(2.75,1.45)(3,1.8)}
        \put(0,-1.5){\qbezier(2.5,1.1)(1.75,1.45)(1,1.8)}
        \put(0,-1.5){\qbezier(2.5,2.9)(2.75,2.35)(3,1.8)}
        \put(0,-1.5){\qbezier(2.5,2.9)(1.75,2.35)(1,1.8)}
        \put(0,-1.5){\qbezier(3,1.8)(2,1.8)(1,1.8)}
        \put(0,-1.5){\qbezier(4,2.1)(3,2.05)(2,2)}
        \put(0,-1.5){\qbezier(4,2.1)(3.25,1.6)(2.5,1.1)}
        \put(0,-1.5){\qbezier(4,2.1)(3.25,2.5)(2.5,2.9)}
         \put(0,-1.5){\qbezier(4,2.1)(3.5,1.95)(3,1.8)}
      \end{picture}
    \end{center}
    \caption{Combination of two integral point sets.}
    \label{fig2}
  \end{figure}

For the construction of plane integral point sets $\Ps$ in semi-general position at first our used method is to combine two point sets consisting of $n-1$ points having $n-2$ points in common to an integral point set consisting of $n$ points, see Figure \ref{fig2}. We remark that this is similar to the approach of \cite{hab_kemnitz}.

As an algorithm for the combination of integral point sets we use a variant of orderly generation  \cite{0412.05006,orderly,mckay,winner}. The big advantage of orderly generation is that the isomorphism problem can be solved without comparing every constructed pair of discrete objects. There is no need to access a large set of constructed structures during the algorithm, so memory is not the bottleneck any more, which is the case for other types of enumeration algorithms. For more information about the concept of orderly generation we refer to \cite{winner} or the more general overview on enumeration algorithms \cite{oestergard}.

Since our variant of orderly generation is applicable for the enumeration of universal discrete structures which can be 
described by an equivalency relation we have to be more general and technical in our description.

At first we have to describe the objects of a discrete structure which we like to construct as an equivalency relation 
$\simeq$ on an ordered list $L$. In our case $L$ is the set of distance matrices which correspond to plane integral 
point sets. To describe $\simeq$ we introduce a mapping $\rho:L\rightarrow\mathbb{N}$. Here, we simply map onto the number $q$ of rows or columns of the distance matrix. This number also equals the number of points of the  corresponding point set. We denote the total ordering of $L$ by $\prec$. If $\rho(l_1)<\rho(l_2)$ or $\rho(l_1)>\rho(l_2)$ for   $l_1,l_2\in L$ we define $l_1\prec  l_2$ or $l_1\succ  l_2$, respectively. In the remaining case $\rho(l_1)=\rho(l_2)$ we continue $\prec$ by a column-lexicographic ordering on the upper right triangle of $l_1$ and $l_2$. By $L_q$ we denote the ordered list containing all elements $l\in L$ with $\rho(l)=q$. If $\rho(l_1)\neq\rho(l_2)$ for $l_1,l_2\in L$ we define $l_1\not\simeq l_2$. Otherwise we define $\simeq$ by the natural group action \cite{0951.05001} of the symmetric group $S_q$ on the $q=\rho(l_1)=\rho(l_2)$ points of $l_1$ and $l_2$.

Because we want to combine only pairs of  distance matrices $l_1,l_2$ consisting of $n-1$ points which have $n-2$ points in common we need a mapping $\pfo$ which deletes the last row and the last column of a distance matrix so that we can rewrite our condition to $\pfo l_1=\,\pfo l_2$. The combination of two distance matrices $l_1$ and $l_2$ itself is done by a mapping $\Gamma_q:L_q\times L_q\rightarrow L_{q+1}^\star$ where $L_{q+1}^\star$ denotes the set of lists of arbitrary length containing elements from $L_{q+1}$. We define $\Gamma_q(l_1,l_2)$ by an example. Let
$l_1=
  \begin{pmatrix}
     0 & 4 & 4 \\[-2mm]
     4 & 0 & 2 \\[-2mm]
    4 & 2 & 0
  \end{pmatrix}
$
and
$l_2=
  \begin{pmatrix}
    0 & 4 & 2 \\[-2mm]
    4 & 0 & 4 \\[-2mm]
    2 & 4 & 0
  \end{pmatrix} 
$
be distance matrices with $\pfo l_1=\,\pfo l_2$ then each distance matrix $\Delta$ in $\Gamma_q(l_1,l_2)$ has the  shape 
$\Delta=
  \begin{pmatrix}
    0 & 4 & 4 & 2 \\[-2mm]
    4 & 0 & 2 & 4 \\[-2mm]
    4 & 2 & 0 & \star \\[-2mm]
    2 & 4 & \star & 0
  \end{pmatrix}
$
where the $\star$ stands for an arbitrary value. In general the distance matrix $\Delta$ is obtained from $l_1$ and $l_2$ by appending the last row and the last column of $l_2$ to $l_1$. The proper values of $\star$ can be calculated by a little computation in the Euclidean metric. It is easy to see that $\star$ can take at most two different values. Additionally we require that the values of $\star$ are integers and that $\Gamma_q(l_1,l_2)$ is ordered by $\prec$. In the above example we have
$$\Gamma_q(l_1,l_2)=\left(\,
    \begin{pmatrix}
      0 & 4 & 4 & 2 \\[-2mm]
      4 & 0 & 2 & 4 \\[-2mm]
      4 & 2 & 0 & 3 \\[-2mm]
      2 & 4 & 3 & 0
    \end{pmatrix}
\,\right).$$
The last ingredient for an orderly generation algorithm is a definition of canonicity. In each equivalence class we
have to mark exactly one element. This element is called canonical. There are several ways to define canonical elements.
We say that an element $l\in L_q$ is canonical if $l\succeq \sigma(l)\,\,\forall\,\sigma\in S_q$, i.e. it is the largest
object in its equivalence class. For our variant of orderly generation we need also another definition. We call an
element $\,l\in L_q\,$ semi-canonical if $\pfo l\succeq \pfo \sigma(l)\,\,\forall \sigma\in S_q$. So each canonical element is also semi-canonical. For algorithmic purposes we define a function $\chi$ which maps an element $l\in L$ into the set $\{canonical,semi\text{\it{-}}canonical,none\}$ where $\chi(l)=semi\text{\it{-}}canonical$ if $l$ is semi-canonical and not canonical. Our aim is to construct complete lists $\mathcal{L}_q$ of the semi-canonical plane integral point sets in semi-general position consisting of $q$ points. We suppose that we have already a list $\mathcal{L}_3$ of all integral semi-canonical triangles with given diameter which can be obtained by a simple double loop. To recursively construct the lists $\mathcal{L}_{q+1}$ we apply the following orderly algorithm:

\begin{algo}{$\,$\\Generation of semi-canonical integral point sets in semi-general position\\}
  \label{algo_orderly}
  \textit{input:} $\mathcal{L}_q$, $\Gamma_q$, $\pfo$ \,, $\prec$\\
  \textit{output:} $\mathcal{L}_{q+1}$\\
  {\bf begin}\\
  \Tsm $\mathcal{L}_{q+1}=\emptyset$\\
  \Tsm{\bf loop over} $l_1\in \mathcal{L}_q$, $\chi(l_1)=canonical$ {\bf do}\\
  \Tsm\Tsm{\bf loop over} $l_2\preceq l_1$, $l_2\in \mathcal{L}_q$, $\pfo l_2=\,\pfo l_1$ {\bf do}\\
  \Tsm\Tsm\Tsm{\bf loop over} $y\in\Gamma_q(l_1,l_2)$ {\bf do}\\
  \Tsm\Tsm\Tsm\Tsm{\bf if} $\chi(y)\neq none$ {\bf then} append $y$ to $\mathcal{L}_{q+1}$ {\bf end}\\
  \Tsm\Tsm\Tsm{\bf end}\\
  \Tsm\Tsm{\bf end}\\
  \Tsm{\bf end}\\ 
  {\bf end}
\end{algo}

The plane integral point set in semi-general position given by the distance matrix
$\Delta=\begin{pmatrix}0&100&89&21\\[-2mm]100&0&21&89\\[-2mm]89&21&0&82\\[-2mm]21&89&82&0\end{pmatrix}$ shows why semi-canonical elements are needed for an exhaustive generation. We notice that $\Delta$ is canonical and can be only combined with a canonical and a semi-canonical triangle. We leave the proof of the correctness of Algorithm \ref{algo_orderly} to the reader and also refer to \cite{phd_kurz}, as it is a bit technical but not difficult.

For the canonicity check $\chi$ we use backtracking with isomorphism pruning in the general case. Because we
have to check $4\times 4$-matrices very often we have developed a fast algorithm for this special case,
see \cite{tetrahedron}. It needs at most $6$ integer comparisons to decide whether a given $4\times 4$-matrix
is canonical, semi-canonical, or none of both. In our special case of integral point sets we can improve Algorithm \ref{algo_orderly} by using the characteristic of an integral point set, see Section \ref{section_characteristic}.

\subsection{Clique search}
\label{sub_sec_clique_search}

In this subsection we will present a hybrid construction algorithm that combines the orderly algorithm of the previous subsection and clique search to search for integral point sets with large diameter. We remark that clique search is a common technique in extremal combinatorics and was also utilized before in the construction of $3$-dimensional integral point sets, see \cite{dipl_piepmeyer}. Suppose we are given an integral triangle $\Delta$. Via orderly generation we can construct all integral point sets $\Ps_i$ in semi-general position consisting of $4$ points with $\pfo \Ps_i=\Delta$. This means that each point set $\Ps_i$ consists of $\Delta$ and a further point $v_i$. The next step is to build up a graph $\mathcal{G}_\Delta$ with $v_i$ as its vertices. We define $\{v_i,v_j\}$ to be an edge in $\mathcal{G}_\Delta$ iff the distance between $v_i$ and $v_j$ is integral and $v_i,v_j$ are not collinear with a point of $\Delta$.

Obviously for every plane integral point set $\Ps$ in semi-general position there exists an integral triangle $\Delta$ such that $\Ps$ corresponds to a clique of $\mathcal{G}_\Delta$. In the other direction every clique of  $\mathcal{G}_\Delta$ corresponds to a plane integral point set but it could happen that three points $v_h$, $v_i$, and $v_j$ are collinear. So we produce only candidates of plane integral point sets which have to be checked whether they are in semi-general position.

Our hybrid algorithm works as follows. For a given diameter $d$ and a lower bound $b$ for the number of points we loop over all integral triangles with diameter $d$. Then we determine the vertices of the graph $\mathcal{G}_\Delta$ by orderly generation and insert the edges. Here, we can use the lower bound $b$ to shrink the graph by deleting edges with at most $b-5$ neighbors. On the resulting graph we perform a clique search  using \textsc{Cliquer} \cite{cliquer,1019.05054} or an implementation of the Bron-Kerbosch algorithm \cite{0261.68018} to generate the cliques of size at least $b-3$. As a last step we reconstruct from the vertices of each clique and the triangle $\Delta$ a
plane integral point set $\Ps$ and check if it is in semi-general position.

\begin{figure}[ht]
  \begin{center}
    \setlength{\unitlength}{0.5cm}
    \begin{picture}(5,5.4)
      \put(0,0.8){\line(1,0){5}}
      \put(0,0.8){\line(3,4){3}}
      \put(5,0.8){\line(-1,2){2}}
      \put(-0.4,-0.1){$A$}
      \put(4.9,-0.1){$B$}
      \put(2.8,5.0){$C$}
      \put(-0.1,0.75){\circle*{0.35}}
      \put(3,4.75){\circle*{0.35}}
      \put(5.1,0.75){\circle*{0.35}}
      \put(0.5,0.75){\circle*{0.35}}
      \put(1.2,0.75){\circle*{0.35}}
      \put(2.1,0.75){\circle*{0.35}}
      \put(2.7,0.75){\circle*{0.35}}
      \put(4,0.75){\circle*{0.35}}
      \put(4.5,0.75){\circle*{0.35}}
    \end{picture}\\[2mm]
  \end{center}
  \caption{Integral points on the side of a triangle.}
  \label{fig_points_on_side_of_triangle}
\end{figure}
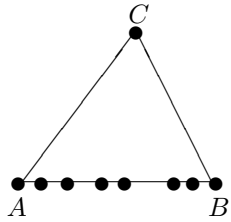

If we want to generate plane integral point sets in arbitrary position we have to modify our approach slightly since it does not produce all possibilities of plane integral point sets consisting of $4$ points which can be seen as follows. Suppose we are given a triangle $\Delta$ as in Figure \ref{fig_points_on_side_of_triangle} where $\overline{AB}$ 
is the largest side. Algorithm \ref{algo_orderly} joins $\Delta$ with all other possible triangles along the side 
$\overline{AB}$ so it cannot generate point sets with further points on the line $l$ through $A$ and $B$.

The situation can be cleared easily. We may simply test all points on $l$ with integral distances to the endpoints of $l$ whether their distance to the third point $C$ of $\Delta$ is also integral. In Section \ref{sec_decomposition_lemma} we will describe a more sophisticated algorithm for this purpose. So we simply add those points on $l$ and the corresponding edges to $\mathcal{G}_\Delta$.

We remark that in our computer calculations the major part of the running time of our hybrid algorithm is used for the orderly generation. This part could be replaced by a more direct algorithm which construct all points $v_i$ which have integral distances to the three corners of a $\Delta$ without the restriction of the diameter, see \cite{paper_axel} for details. We do not use this algorithm for our purpose because it runs more slowly.

\section{Characteristic of a plane integral point sets}
\label{section_characteristic}

% If we use Algorithm \ref{algo_orderly} for the construction of plane integral point sets consisting of $4$ points in
% semi-general position with diameter at most $d$ then we have to expect a running time of $O(d^5)$ because there are
% $O(d^2)$ integral triangles with diameter $d$. 
Due to the Heron formula
$$
  A_\Delta(a,b,c)=\frac{\sqrt{(a+b+c)(a+b-c)(a-b+c)(-a+b+c)}}{4}
$$
for the area $A_\Delta$ of a triangle with side lengths $a$, $b$, $c$ we can write the area of a non-degenerate   triangle with integral side lengths uniquely as $A=q\sqrt{k}$ with a rational number $q$ and a squarefree integer $k$.
The number $k$ is called the characteristic $\car(\Delta)$ of the triangle. The next theorem allows us to  talk about the characteristic $\car(\Ps)$ of a point set $\Ps$.

\begin{theorem}
  \label{car}
  Each non-degenerate triangle $\Delta$ in a plane integral point $\Ps$ set has the same characteristic
  $\car(\Ps)=\car(\Delta)$.
\end{theorem}
\begin{proof}
  See \cite{hab_kemnitz} or \cite{phd_kurz,characteristic}
\end{proof}

Clearly, we modify Algorithm \ref{algo_orderly} and combine only integral point sets with equal characteristic. To be able to measure the complexity of the modified algorithm and for some number theoretical insight we introduce the function $\psi(d,k)=$
$$
  \left|\left\lbrace (a,b)\,\Big|\,
  \begin{array}{l}
     a,b\in\{1,2,\dots, d\},\, a+b>d,
    \exists w\in\mathbb{N}:kw^2=\\(a+b+d)(a+b-d)(a-b+d)(-a+b+d)
  \end{array} 
  \right\rbrace\right|
$$
which counts the number of integral semi-canonical triangles with diameter $d$ and characteristic $k$ if we set $\psi(d,k)=0$ for non-square\-free numbers $k$. We remark that the the running time of the modified orderly algorithm for the generation of integral point sets consisting of $4$ points is given by $\sum_k{\psi(d,k)+1\choose 2}$.

\begin{figure}[ht]
  \centering
  \includegraphics[width=10cm]{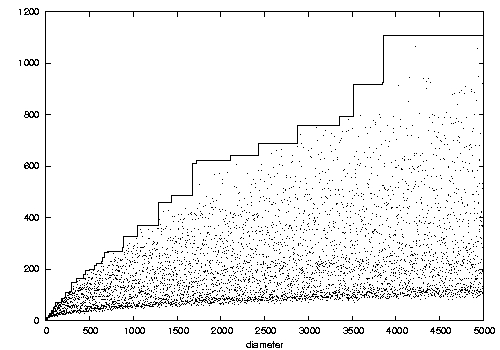}
  \caption{Occurring values $\psi(d,k)$ for $d\le 5000$.}
  \label{fig_psi_tilde}
\end{figure}

In Figure \ref{fig_psi_tilde} we have depicted the occurring values of $\psi(d,k)$ for $d\le 5000$. We introduce the  function $\tilde{\psi}(d)=\max_k\{\psi(d,k)\mid k\in\mathbb{N}\}$. 
%together with the estimation 
%$$
%  \Psi(d)\le\frac{\tilde{\psi}(d)(\tilde{\psi}(d)+1)}{2}\frac{d(d+1)}{2\tilde{\psi}(d)}=
%  \frac{d(d+1)(\tilde{\psi}(d)+1)}{4},
%$$
%where we have used the fact that there are $\frac{d(d+1)}{2}$ semi-canonical integral triangles with
%diameter $d$.
A parameter solution for the set of all integral triangles with characteristic
$k$ given in \cite{hab_kemnitz} can be rewritten to:
\begin{theorem}
  \label{thm_parameterization}
  For all integral triangles with side lengths $a$, $b$, $c$ and characteristic $k$ there exists at least 
  one integer tuple $(p,q,h,i,j)$ fulfilling\\[-12mm]
  \begin{eqnarray*}
    a&=&\frac{ph(i^2+kj^2)}{q},\\
    b&=&\frac{pi(h^2+kj^2)}{q},\\
    c&=&\frac{(i+h)(ih-kj^2)}{q},
  \end{eqnarray*}\\[-5mm]
  $\gcd(p,q)=\gcd(h,i,j)=1$, $i\ge h$, and $ih>kj^2$.
\end{theorem}

Using Theorem \ref{thm_parameterization} and some cumbersome technical computation one can deduce
$\tilde{\psi}(d)\in O\!\left(d^{1+\frac{c}{\log\log d}}\right)$ for a suitable constant $c$. For the details we refer the interested reader to \cite{phd_kurz}, where also more detailed considerations on the running time of Algorithm \ref{algo_orderly} and some more numerical data can be found.

%For the details of the computation we refer to \cite{phd_kurz}.
%For some numerical data see Figure \ref{fig_psi_tilde}. By combining our two estimations we receive 
%
%  \begin{lemma}
%    $$\Psi(d)\in O\!\left(d^{3+\frac{c}{\log\log d}}\right).$$
%  \end{lemma}
%
%
% Thus our modified orderly algorithm has provable running time $O\!\left(d^{4+\frac{c}{\log\log d}}\right)$ and
% conjectured running time $O\!\left(d^{3+\frac{c}{\log\log d}}\right)$ for the generation of integral point sets in
% semi-general position consisting of $4$ points with diameter at most $d$.

% In contrast to Lemma \ref{lemma_tilde_psi} the following lemma states that there exist many different characteristics
% for integral triangles with diameter $d$.
%

In the other direction we have:\\[-8mm]
\begin{lemma}
  The number of different characteristics of integral triangles with diameter $d$ is in $O(d^2)$ and in
  $\Omega\!\left(\frac{d^2}{(\log d)^2}\right)$.
\end{lemma}
\begin{proof}
  Since there are $O(d^2)$ integral triangles with diameter $d$ we have the stated trivial upper bound.
  %The $\lfloor\frac{(d+1)^2}{4}\rfloor$ non-isomorphic integral triangles with diameter $d$ are clearly an
  %upper bound for the number of different characteristics. 
  For the other direction we choose for suitable large $d$ two
  primes $p_1,p_2$ fullfilling $\frac{9}{4}d< p_1<\frac{10}{4}d$ and $\frac{5}{4}d<p_2<\frac{6}{4}d$. With this we set
  $a=d$, $b=\frac{p_1+p_2}{2}-d$, and $c=\frac{p_1-p_2}{2}$. Thus $$a+b+c=p_1,\,\, a+b-c=p_2,\,\, \frac{3}{4}d<b<d,\, \,
  \mbox{ and }\,\,\frac{3}{8}d<c<\frac{5}{8}d.$$
  Because $p_1$ and $p_2$ have to be odd for big enough $d$ and because $b+c=p_1-d>d=a$ the values $a>b>c$ are integers
  and fulfill the triangle conditions. Since $\frac{1}{2}d<a-b+c<\frac{3}{4}d$ and $\frac{1}{4}d<-a+b+c<\frac{1}{2}d$
  the characteristic of the triangle with side lengths $a$, $b$, and $c$ is divisible by $p_1p_2$. Due to the prime
  number theorem we have $\Omega\!\left(\frac{d}{\log d}\right)$ choices for $p_1$ and $p_2$ each. Thus there are at
  least $\Omega\!\left(\frac{d^2}{(\log d)^2}\right)$ different characteristics.
\end{proof}

% We will come back to the properties of the characteristic of plane integral point sets later on when we will propose
% a strategy to prove a good lower bound for the minimum diameter $d(2,n)$. As a last remark 
We would like to mention that the first author has recently generalized the definition of the characteristic of a plane integral point set to $\di$-dimensional integral point sets and has proven an analog theorem to Theorem \ref{car}, see \cite{phd_kurz,characteristic}.

%----------------------------------------------------------------------------------------------------------------------
\section{Minimum diameter of plane integral point sets}
\label{sec_minimum diameter}

With the algorithms of Section \ref{sec_exhaustive_search} at hand we are now state the results of our exhaustive search for plane integral point sets with minimum diameter.

\subsection{Plane integral point sets in general position}
\label{subsection_general_position}
To construct plane integral point sets in general position we need a check to decide whether three points lie on a line or four   points are located on a circle. For the first check we can use the triangle inequalities in the degenerated case and for the second check we can use Ptolemy's theorem.

We have implemented the algorithm described in Subsection \ref{sub_sec_orderly} and our computers have constructed all
plane integral point sets in general position with diameter at most $30\,000$. We have only found one such point set consisting of seven points, which proves Theorem \ref{thm_general_position}. For details, pictures, and a second example of diameter $66\,810$ we refer to \cite{kreisel}.

\subsection{Plane integral point sets in semi-general position}
\label{subsection_semi_general_position}
We have also applied our orderly algorithm on the construction of plane integral point sets in semi-general position.
Our available computer power has allowed us to generate all such point sets with diameter at most $5\,000$ and thus to
determine $\overline{d}(2,n)$ up to $n=24$. 
% The reason why our computers were not able to construct all point sets in semi-general position up to diameter
% $20\,000$ is due to the fact that as $n$ grows a plane integral point set consisting  of $n$ points contains many
% non-isomorphic integral point sets with the same diameter consisting of $\lceil\frac{d}{2}\rceil$ 
% points. We will be more explicit in the case of plane integral point sets in arbitrary position later on.
%
With the hybrid construction algorithm of Subsection \ref{sub_sec_clique_search} we were able to enumerate all plane integral point sets in semi-general position with minimum diameter $\overline{d}(2,n)$ up to $\diam(\Ps)=20\,000$ leading to
\begin{eqnarray*}
  \!\!\!\!\!\!\left(\overline{d}(2,n)\right)_{n=10,\dots,36}=105,105,105,532,532,735,735,735,735,\\
  1995,1995,1995,1995,1995,1995,9555,9555,9555,10672,\\
  13975,13975,13975,13975,13975,13975,13975,13975.
\end{eqnarray*}

We remark that our hybrid algorithm of Subsection \ref{sub_sec_clique_search} did not produce any candidates with collinear triples.

\begin{table}[ht]
  \begin{center}
    \begin{tabular}{|r|r|r|r||r|r|r|r|}
      \hline
      $z$ & $k$ & $|\Ps|$ & $\diam(\Ps)$ & $z$ & $k$ & $|\Ps|$ & $\diam(\Ps)$ \\
      \hline
      $1$ & $3$ &   3 &          1 &
      %$2^4\cdot 17\cdot 19\cdot 23$ 
      118864 & $15$ &  40 &     61\,375 \\
      %$2^3$
      8& $15$ &   4 &          4 &
      %$7\cdot 13\cdot 19\cdot 31$
      53599 & $3$ &  48 &     61\,880 \\
      $7$ & $3$ &   6 &          8 &
      %$7^2\cdot 13^2\cdot 19$
      157339 & $3$ &  54 &    181\,675 \\
      %$2^6$ 
      64 & $15$ &   7 &         33 &
      %$2^6\cdot 17\cdot 19\cdot 23$
      475456 & $15$ &  56 &    245\,518 \\
      %$7^2$
      49 & $3$ &   9 &         56 &
      %$7^2\cdot 13\cdot 19\cdot 31$
      375193 & $3$ &  72 &    433\,225 \\
      %$7\cdot 13$ 
      91 & $3$ &  12 &        105 &
      %$2^4\cdot 17\cdot 19\cdot 23\cdot 31$
      3684784 & $15$ &  80 &  1\,902\,813 \\
      %$2^6\cdot 11$
      704 & $7$ &  14 &        532 &
      %$7\cdot 13\cdot 19\cdot 31\cdot 37$
      1983163 & $3$ &  96 &  2\,289\,957 \\
      %$7^2\cdot 13$ 
      637 & $3$ &  18 &        735 &
      %$7^2\cdot 13^2\cdot 19\cdot 31$
      4877509 & $3$ & 108 &  5\,632\,056 \\
      %$7\cdot 13\cdot 19$ 
      1729 & $3$ &  24 &      1\,995 &
      %$2^6\cdot 17\cdot 19\cdot 23\cdot 31$
      14739136 & $15$ & 112 &  7\,611\,252 \\
      %$7^2\cdot 13^2$ 
      8281 & $3$ &  27 &      9\,555 &
      %$7^2\cdot 13\cdot 19\cdot 31\cdot 37$
      13882141 & $3$ & 144 & 16\,029\,704 \\
      %$2^6\cdot 17\cdot 19$
      20672 & $15$ &  28 &    10\,672 &
      %$7\cdot 13\cdot 19\cdot 31\cdot 37\cdot 43$
      85276009 & $3$ & 192 & 98\,468\,151 \\
      %$7^2\cdot 13\cdot 19$
      12103 & $3$ &  36 &     13\,975 &
      &     &     &            \\
      \hline
    \end{tabular}
    \caption{Upper bounds for $\overline{d}(2,n)$ from point sets on circles with radius $r=\frac{z}{k}\sqrt{k}$.}
    \label{table_upper_bounds_semi-general_2}
  \end{center}
\end{table}

%\begin{conjecture}
%  \label{conjecture_n3ol}
%  The points of a plane integral point set in semi-general position with minimum diameter are located on a circle.
%\end{conjecture}

With the minimal examples up to to $n=36$ at hand we could check Theorem \ref{thm_n3ol} which states that the points of 
integral point sets in semi-general position with minimum diameter are located on circles with radius 
$r=\frac{z}{k}\sqrt{k}$ where $k\in\{3,7,15\}$ and $z$ is an integer with many prime factors. With this we are motivated to conjecture this pattern in general, see Conjecture \ref{conjecture_on_line}. Searching for these special types of integral point sets yields the upper bounds given in Table \ref{table_upper_bounds_semi-general_2}.

\begin{conjecture}
  The upper bounds from Table \ref{table_upper_bounds_semi-general_2} are sharp.
\end{conjecture}

We remark that the authors of  \cite{minimum_diameter} describe a construction, for $k=3$ and $z=\prod_{i=1}^{r}p_i^{v_i}$ with $p_i\equiv 1\mod 6$ being primes, directly yielding the integral point
sets of Table \ref{table_upper_bounds_semi-general_2} in these cases. Their construction relies on calculations
over the ring $\mathbb{Z}\left[\frac{-1+\sqrt{-3}}{2}\right]$. It is possible that there are similar constructions
for squarefree $k$ or at least $k=7,15$. 

\subsection{Plane integral point sets in arbitrary position}
\label{sub_sec_arbitrary}

The orderly generation by combination approach is limited by the number of substructures as $n$ grows. For example we consider a plane integral point set $\Ps$ consisting of $89$ points where $88$ points are collinear. We will determine $d(2,89)$ and it will turn out that the corresponding point set has the shape of Figure \ref{fig_points_on_side_of_triangle}. The convex hull of $\Ps$ is formed by three points $A$, $B$, and $C$. Let
us assume that the $86$ other points of $\Ps$ are located on the line through $A$ and $B$. Now we consider all point   sets consisting of $A$, $B$, $C$, and $43$ of the other points. Then we have the maximum number of ${86 \choose 43}\approx 6.6\cdot 10^{24}$ possibilities. Because orderly generation has to generate all these point sets to finish in $\Ps$ it is beyond our means to determine $d(2,89)$ with this method.

%\begin{figure}[ht]
%  \begin{center}
%    \setlength{\unitlength}{0.5cm}
%    \begin{picture}(5,5.4)
%      \put(0,0.8){\line(1,0){5}}
%      \put(0,0.8){\line(3,4){3}}
%      \put(5,0.8){\line(-1,2){2}}
%      \put(-0.4,-0.1){$A$}
%      \put(4.9,-0.1){$B$}
%      \put(2.8,5.0){$C$}
%      \put(-0.1,0.75){\circle*{0.35}}
%      \put(3,4.75){\circle*{0.35}}
%      \put(5.1,0.75){\circle*{0.35}}
%      \put(0.5,0.75){\circle*{0.35}}
%      \put(1.2,0.75){\circle*{0.35}}
%      \put(2.1,0.75){\circle*{0.35}}
%      \put(2.7,0.75){\circle*{0.35}}
%      \put(4,0.75){\circle*{0.35}}
%      \put(4.5,0.75){\circle*{0.35}}
%    \end{picture}\\[2mm]
%  \end{center}
%  \caption{Integral points on the side of a triangle.}
%  \label{fig_points_on_side_of_triangle}
%\end{figure}

So again we had to use the hybrid algorithm of Subsection \ref{sub_sec_clique_search}. By doing an exhaustive search up to diameter $10\,000$ we were able to determine the minimum diameter $d(2,n)$ up to $n=122$ points.
% The computational improvement compared to \cite{integral_distances_in_point_sets} amounts to a factor of
% approximately $\left(\frac{10\,000}{29}\right)^3>40\,000\,000$. 
Checking the minimal examples uncovers that for $9\le n\le 122$ points in each case $n-1$ 
points are collinear, which proves Theorem \ref{thm_on_line} and motivates Conjecture \ref{conjecture_on_line}.

%\begin{conjecture}
%  \label{conjecture_on_line}
%  There exists a constant $\kappa$ for $n\in\mathbb{N}$ so that each plane integral point set with minimum diameter 
%  consisting of $n$ points contains a subset of $n-\kappa$ collinear points.
% \end{conjecture}
%
%----------------------------------------------------------------------------------------------------------------------  

\section{Plane integral point sets with many points on a line}
\label{sec_decomposition_lemma}
We have seen in the last subsection that plane integral point sets with many collinear points are interesting objects to study because plane integral point sets with minimum diameter seem to belong to this class. At first we present an important link between plane integral point sets with $n-1$ collinear points and the decompositions of a certain integer, the decomposition number, into two factors.

\begin{definition}
  The decomposition number $D$ of an integral triangle with side lengths $a$, $b$, and $c$ is given by
  $$
    D=\frac{(a+b+c)(a+b-c)(a-b+c)(-a+b+c)}{\gcd(b^2-c^2+a^2,2a)^2}.
  $$
\end{definition}

\begin{figure}[ht]
  \begin{center}
    \setlength{\unitlength}{0.7cm}
    \begin{picture}(9,4.5)
      \put(0,0){\line(1,0){9}}
      \put(0,0){\line(5,4){5}} 
      \put(1,0){\line(1,1){4}}
      \put(2,0){\line(3,4){3}}
      \put(4,0){\line(1,4){1}}
      \put(6,0){\line(-1,4){1}}
      \put(7,0){\line(-1,2){2}}
      \put(9,0){\line(-1,1){4}}
      \put(0,0){\circle*{0.3}} 
      \put(1,0){\circle*{0.3}} 
      \put(2,0){\circle*{0.3}} 
      \put(4,0){\circle*{0.3}} 
      \put(6,0){\circle*{0.3}} 
      \put(7,0){\circle*{0.3}}  
      \put(9,0){\circle*{0.3}} 
      \put(5,4){\circle*{0.3}} 
      \put(5,0){\line(0,1){4}}
      \put(0.1,-0.5){$a_3$}
      \put(1.1,-0.5){$a_2$}
      \put(2.8,-0.5){$a_1$}
      \put(4.4,0.2){$q$}
      \put(5.2,0.2){$q'$}
      \put(5.1,1.4){$h$}
      \put(4.7,-0.5){$a_0$}
      \put(6.2,-0.5){$a_1'$}
      \put(7.6,-0.5){$a_2'$}
      \put(1.6,2){$b_3$}
      \put(0.8,0.4){$b_2$}
      \put(1.7,0.4){$b_1$}
      \put(3.5,0.4){$b_0$}
      \put(6.05,0.4){$b_0'$}
      \put(7.0,0.4){$b_1'$}
      \put(7,2){$b_2'$}
    \end{picture}\\[2mm]
  \end{center}
  \caption{Plane integral point set $\mathcal{P}$ with $n-1$ points on a line.}
  \label{fig_on_a_line}
\end{figure}
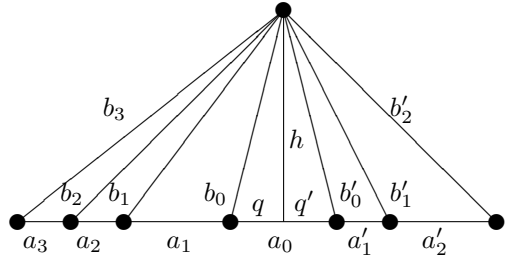

\begin{lemma}{\bf(Decomposition Lemma)\\}
  \label{lemma_decomposition}
  \noindent
  The distances of a plane integral point set $\mathcal{P}$ consisting of $n$ points where a subset of $n-1$ points is 
  collinear correspond to decompositions of the decomposition number $D$ of the largest triangle of $\mathcal{P}$ 
  into two factors.
\end{lemma}
\begin{proof}
  We use the notation of Figure \ref{fig_on_a_line}, let $s$ be the largest index of $a_i$, let $t$ be the largest
  index of $a'_i$, and set $$c_i=q+\sum_{j=1}^ia_j\quad\mbox{for}\quad 0\le i\le s,\,
  c_i'=q'+\sum_{j=1}^ia_j'\quad\mbox{for}\quad 0\le i\le t.$$
  %Pythagoras' Theorem yields $c_{i+1}^2+h^2=b_{i+1}^2$ and $c_i^2+h^2=b_i^2$ for $0\le i< s$. We subtract these
  %equations from each other and get $$a_{i+1}^2+2a_{i+1}\sum_{j=1}^ia_j+2a_{i+1}q=b_{i+1}^2-b_i^2\,.$$
  %Because the $a_i$ and the $b_i$ are  positive integers we have $2a_{i+1}q\in\mathbb{N}$ for $0\le i<s$, 
  %i.e. $2a_sq\in\mathbb{N}$. % and therefore $2\gcd(a_1,a_2,\dots,a_s)q\in\mathbb{N}$.
  Let us associate $b$ with $b_s$, $c$ with $b_t'$, and $a$ with $\sum_{i=1}^sa_i+a_0+\sum_{i=1}^ta_i'$. With this 
  we have $$c_s=q+\sum\limits_{j=1}^sa_j=\frac{b^2-c^2+a^2}{2a}.$$
  By defining $g:=\frac{2a}{\gcd(b^2-c^2+a^2,2a)}$ we obtain $gc_s\in\mathbb{N}$ and $gq\in\mathbb{N}$. From
  $q+q'=a_0\,\,\in\mathbb{N}$ we conclude $gq'\in\mathbb{N}$. Thus $gc_i,gc'_j\in\mathbb{N}$ for all possible 
  indices. An application of Pythagoras' Theorem yields $c_{i+1}^2+h^2=b_{i+1}^2$ and $c_i^2+h^2=b_i^2$ for
  $0\le i< s$. We conclude $$g^2h^2=(gb_i+gc_i)(gb_i-gc_i)=(gb_j'+gc_j')(gb_j'-gc_j').$$
  Since the $b_i$, $b'_j$ are integers we can obtain the possible values for $c_i$ and $c_i'$ by decomposing 
  $g^2h^2$ into two factors.

  Due to the Heron formula $16A_\Delta^2=(a+b+c)(a+b-c)(a-b+c)(-a+b+c)$ and the formula for the area of a triangle
  $2A_\Delta=ah$ we finally get
  \begin{eqnarray*}
    g^2h^2=\frac{g^2(a+b+c)(a+b-c)(a-b+c)(-a+b+c)}{4a^2}=\\
    =\frac{(a+b+c)(a+b-c)(a-b+c)(-a+b+c)}{\gcd(b^2-c^2+a^2,2a)^2}=D\,.
  \end{eqnarray*}
\end{proof}

With the aid of the Decomposition Lemma and Theorem \ref{thm_on_line} we are able to describe the plane integral points with minimum diameter for $9\le n\le 122$ points in a very compact manner by giving $n$, $D=g^2h^2$, and $g$. For $9\le n\le 20$ see Table \ref{table_minimum_planar_1} and for the complete listing see \cite{own_hp2}. We remark that for $n\le 122$ we only have $g=2$ for $n=9,10,11$ otherwise we have $g=1$. The Decomposition Lemma can also be used as a heuristic to determine good upper bounds for $d(2,n)$. One only has to loop over Decompositions numbers $D$ with many divisors and apply Lemma \ref{lemma_decomposition} for the construction of an integral point set. For $n\le 148$ the results of this heuristic are given in \cite{own_hp2}.

\begin{table}[ht]
  \begin{center}
    \begin{tabular}{|r|r|c|r||r|r|c|r||r|r|c|}
      \hline
      $n$ & $\!d(2,n)\!$ & $D=g^2h^2$ & $g$ & $n$ & $\!d(2,n)\!$ & $D=g^2h^2$ &$g$\\%& $n$ & $\!d(2,n)\!$ & $D=g^2h^2$ \\
      \hline
       9 &   29 & $3^2\cdot 5\cdot 7$ & 2 &
      15 &  104 & $2^5\cdot 3^2\cdot 5$ & 1\\
      %21 &  212 & $2^6\cdot 3\cdot 5\cdot 7$ \\
      10 &   40 & $3^2\cdot 5\cdot 7$  & 2 &
      16 &  121 & $2^5\cdot 3^2\cdot 5$ & 1\\
      %21 &  212 & $2^3\cdot 3^2\cdot 5\cdot 7$ \\
      11 &   51 & $3^2\cdot 5\cdot 7$ & 2 &
      17 &  134 & $2^5\cdot 3^2\cdot 5$ & 1\\
      %22 &  228 & $2^3\cdot 3^2\cdot 5\cdot 7$ \\
      12 &   63 & $2^5\cdot 3\cdot 5$ & 1 &
      18 &  153 & $2^5\cdot 3^2\cdot 5$ & 1\\
      %23 &  244 & $2^3\cdot 3^2\cdot 5\cdot 7$ \\
      13 &   74 & $2^5\cdot 3\cdot 5$ & 1 &
      18 &  153 & $2^6\cdot 3^2\cdot 5$ & 1\\
      %24 &  272 & $2^3\cdot 3^2\cdot 5\cdot 7$ \\
      13 &   74 & $2^5\cdot 3^2\cdot 5$ & 1 &
      19 &  164 & $2^6\cdot 3^2\cdot 5$ & 1\\
      %25 &  288 & $2^3\cdot 3^2\cdot 5\cdot 7$ \\
      14 &   91 & $2^5\cdot 3^2\cdot 5$ & 1 &
      20 &  196 & $2^6\cdot 3^2\cdot 5$ & 1\\
      %26 &  319 & $2^3\cdot 3^2\cdot 5\cdot 7$ \\
      \hline 
    \end{tabular}
    \caption{Parameters for plane integral point sets with minimum diameter and $9\le n\le 20$.}
    \label{table_minimum_planar_1}
  \end{center}
\end{table}

We can also utilize the Decomposition Lemma for a fast algorithm to  determine the integral points on a side of an integral triangle $\Delta$. Suppose we are given a triangle with side lengths $a,b,c$ and diameter $d$. We also assume
that we have tabulated the prime factorizations of all integers smaller equal $3d$ in a pre-calculation. So we can
use the formula in the proof of Lemma \ref{lemma_decomposition} to obtain the prime factorization of $g^2h^2$ and loop
through the divisors and determine the suitable points on a side of $\Delta$ in time $O\!\left(d^{\frac{c'}{\log \log d}}\right)$ where $c'$ is a suitable constant.

The last application of the Decomposition Lemma is the proof of Theorem \ref{thm_lower_bound}. Therefore we need the following theorem from number theory:

\begin{theorem}{(\textbf{Theorem 317 \cite{number_theory_hw}})}
  \label{lemma_num_theo}
  There exists a function $m_0(\varepsilon)$ such that for each $\varepsilon>0$ and all $m>m_0(\varepsilon)$
  we have
  $$\tau(m)<2^{(1+\varepsilon)\frac{\log m}{\log\log m}}$$
  where $\tau(m)$ denotes the number of divisors of $m$.
\end{theorem}

  %\begin{theorem}
  %  \label{thm_lower_bound}
  %  For a plane integral point set $\mathcal{P}$ with at least $n^\delta$ collinear points we have
  %  $$\mbox{diam}(\mathcal{P})\ge n^{\frac{\delta}{4\log 2(1+\varepsilon)}  \log\log n}$$ 
  %  where $\varepsilon>0$ and $n\ge n_0(\varepsilon)$.
  %\end{theorem}
%\begin{proof}{(Theorem \ref{thm_lower_bound})}\\
\textbf{Proof of Theorem \ref{thm_lower_bound}.}\\
  Due to the Decomposition lemma and $\max(s+1,t+1)\ge\frac{n^\delta}{2}$ we need at least $\frac{n^\delta}{2}$
  different factorizations of $g^2h^2$ into two factors, and so we have the condition
  $\tau(g^2h^2)\ge \frac{n^\delta}{2}$. With $g^2h^2\le 4\,\mbox{diam}(\mathcal{P})^4$ and Theorem
  \ref{lemma_num_theo} we conclude
  $$
    2^{\frac{(1+\varepsilon')\log( 4\cdot\mbox{\begin{tiny}diam\end{tiny}}(\Ps)^4)}{\log\log( 4\cdot
    \mbox{\begin{tiny}diam\end{tiny}}(\Ps)^4)}}\ge\frac{n^\delta}{2}
  $$
  for $n\ge n_0'(\varepsilon')$ and $\varepsilon'>0$. Because $\mbox{diam}(\Ps)\ge\mbox{d}(2,n)\ge c_1 n$ 
  and with a change from $\varepsilon'$ to $\epsilon>0$ we have
  $$
    2^{\frac{(1+\varepsilon)4\log\mbox{\begin{tiny}diam\end{tiny}}(\Ps)}{\log\log n}}\ge n^\delta
  $$
  for $n\ge n_0(\varepsilon)$, and the theorem follows by an elementary transformation.
  \hfill{$\square$}

\section{Some remarks on lower bounds for $\mathbf{d(2,n)}$}
\label{sec_some_remarks}
 So far we have derived a good lower bound (compared to the upper bound) only in the case when the point set contains
\textit{many} collinear points. It would be nice to get rid of this last condition. We would like to mention a possible
reason why plane integral point sets with minimum diameter seem to have subsets of many collinear points. Theorem \ref{car} states that each non-degenerated triangle of a plane integral point set $\Ps$ has the same characteristic. In Section \ref{section_characteristic} we have observed that integral triangles with equal characteristic and equal diameter are somewhat rare. So a plane integral point set with small diameter is forced to either use isomorphic triangles several times or to have many degenerated triangles. For the second possibility we have the following analysis. The maximum number of non-degenerate triangles is achieved by a point set in semi-general position and consists of ${n\choose 3}=\frac{n(n-1)(n-2)}{6}$ non-degenerate triangles where $n$ is the number of points. The minimum number ${n-1\choose 2}=\frac{(n-1)(n-2)}{2}$ is attained by a point set where $n-1$ points are collinear. The first possibility is ruled out by the following lemma.
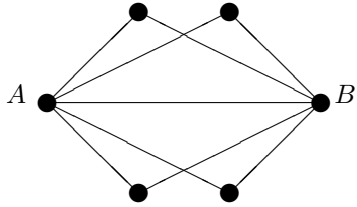
\begin{figure}[ht]
  \begin{center}
    \setlength{\unitlength}{1.2cm}
    \begin{picture}(3,2)
      \put(0,1){\line(1,0){3}}
      \put(0,1){\line(1,1){1}}
      \put(0,1){\line(1,-1){1}}
      \put(0,1){\line(2,1){2}}
      \put(0,1){\line(2,-1){2}}
      \put(3,1){\line(-1,-1){1}}
      \put(3,1){\line(-1,1){1}}
      \put(3,1){\line(-2,-1){2}}
      \put(3,1){\line(-2,1){2}}
      \put(-0.45,1){$A$}
      \put(3.15,1){$B$}
      \put(0,1){\circle*{0.2}}
      \put(3,1){\circle*{0.2}}
      \put(1,2){\circle*{0.2}}
      \put(1,0){\circle*{0.2}}
      \put(2,2){\circle*{0.2}}
      \put(2,0){\circle*{0.2}}
    \end{picture}
  \end{center}
  \caption{Equivalent triangles sharing a common side.}
  \label{fig5}
\end{figure}

\begin{lemma}
  \label{lemma_many_triangles}
  Each plane integral point set consisting of $n$ points with at most $\frac{n}{2}$ points on a line 
  contains a set of at least $\frac{n}{8}$ different integral triangles with equal diameter and equal characteristic.
\end{lemma}
\begin{proof}
  We can choose an arbitrary pair $(A,B)$ of points out of the $n$ points. Because at most
  $\frac{n}{2}$ points lie on a line there are at least $\frac{n}{2}$ points $C$ not on the line through $A$ 
  and $B$. The proof is finished by the fact that at most 4 equivalent triangles can share the side between $A$ and
  $B$, see Figure \ref{fig5}.
\end{proof}

With $\tilde{\psi}(d)\in O\!\left(d^{1+\frac{c}{\log\log d}}\right)$, Lemma \ref{lemma_many_triangles}, and Theorem \ref{thm_lower_bound} one could conclude
$$
   d(2,n)\ge cd^{1-\varepsilon}
$$
for arbitrary $\varepsilon>0$ and a suitable constant $c$, see \cite{phd_kurz} for the details. Unfortunately the bound is slightly less than the bound  of Solymosi \cite{note_on_integral_distances}. So we propose a similar strategy along the same lines. Consider a fix integral triangle $\Delta$ with diameter $d$ and define $\Upsilon(\Delta,d)$ as the number
of canonical integral point sets $\Ps$ in semi-general position with $\pfo \Ps=\Delta$. By $\overline{\Upsilon}(d)$ 
we denote the maximum of $\Upsilon(\Delta,d)$ over all integral triangles $\Delta$  with diameter $d$. Similar to Lemma \ref{lemma_many_triangles} we can prove that a plane integral point set consisting of $n$ points with \textit{few} collinear points and diameter $d$ contains $\Omega(n)$ point sets of the type counted by $\Upsilon(\Delta,d)$. So we have $\frac{n}{c'}\ge \overline{\Upsilon}^{\,-1}(d)$ for a suitable constant $c'$. There is some numerical evidence that
$\overline{\Upsilon}(d)\in O\!\left(d^{\frac{c}{\log\log d}}\right)$ for a suitable constant $c$. If that could be proven we would have a lower bound of $d(2,n)\ge \tilde{c}\cdot n^{c'\log\log n}$ for suitable constants $\tilde{c},c'$, see \cite{phd_kurz} for the details.

%In Figure \ref{fig_upsilon_n3ol} we have plotted some numerical data about 
%  $\overline{\Upsilon}(d)$ which leads us to:

\section{Conclusion}
We have presented some new exact values for the minimum diameters $d(2,n)$, $\overline{d}(2,n)$ and for $\dot{d}(2,7)$ obtained by exhaustive search with custom-made algorithms. Having these new values and the corresponding integral point sets at hand we may speculate about a structure theorem for integral point sets with minimum diameter in arbitrary or semi-general position. We have formulated this as Conjecture \ref{conjecture_on_line} and Conjecture \ref{conjecture_on_circle}, respectively. In Lemma \ref{lemma_decomposition} we have presented an important link
between decompositions of a certain number into two factors and the distances of a plane integral point set consisting
of $n$ points with a subset of $n-1$ collinear points. It also seems that the minimum diameter $\overline{d}(2,n)$ is
dominated by decompositions of certain numbers. Here, \cite{minimum_diameter} gives a first glance at the possibly
underlying rich number theoretic structure, but more research has to be done. The derivation of tight bounds for the
minimum diameter $d(2,n)$ is a challenging task for the future. Our contribution was the resolution of the special case
of point sets with many collinear points.
\bibliography{on_the_minimum_diameter_ars}
\bibdata{on_the_minimum_diameter_ars}
\bibliographystyle{plain}
\end{document}